\documentclass[12pt]{article}
\usepackage{latexsym}
\usepackage{amsmath,amssymb}
\usepackage[dvips]{rotating}

\usepackage[dvips]{graphics,epsfig}

\renewcommand{\labelenumi}{(\roman{enumi})}

\usepackage[thmmarks,amsmath]{ntheorem}

\theorempreskipamount 2ex \theorempostskipamount 3ex \theoremstyle{change}
\theoremsymbol{$\Box$} \theoremseparator{\quad}

\def\rueck{\noindent\hangafter=1 \hangindent=\parindent}

\def\N{{\mathbb{N}}}
\def\R{{\mathbb{R}}}
\def\Z{{\mathbb{Z}}}

\def\eps{{\varepsilon}}

\title{\vspace*{-1cm} Some aspects of extreme value statistics under serial dependence}
\author{Holger Drees\footnote{Dpt.\ of Mathematics, University of Hamburg, Bundesstr.\ 55, 20146 Hamburg, Germany; e-mail:
 {drees@math.uni-hamburg.de}}\\
  University of Hamburg}
\date{}

\begin{document}

\vspace*{-2.5cm}

\nopagebreak

\thispagestyle{empty} \centerline{\thispagestyle{empty}
\includegraphics{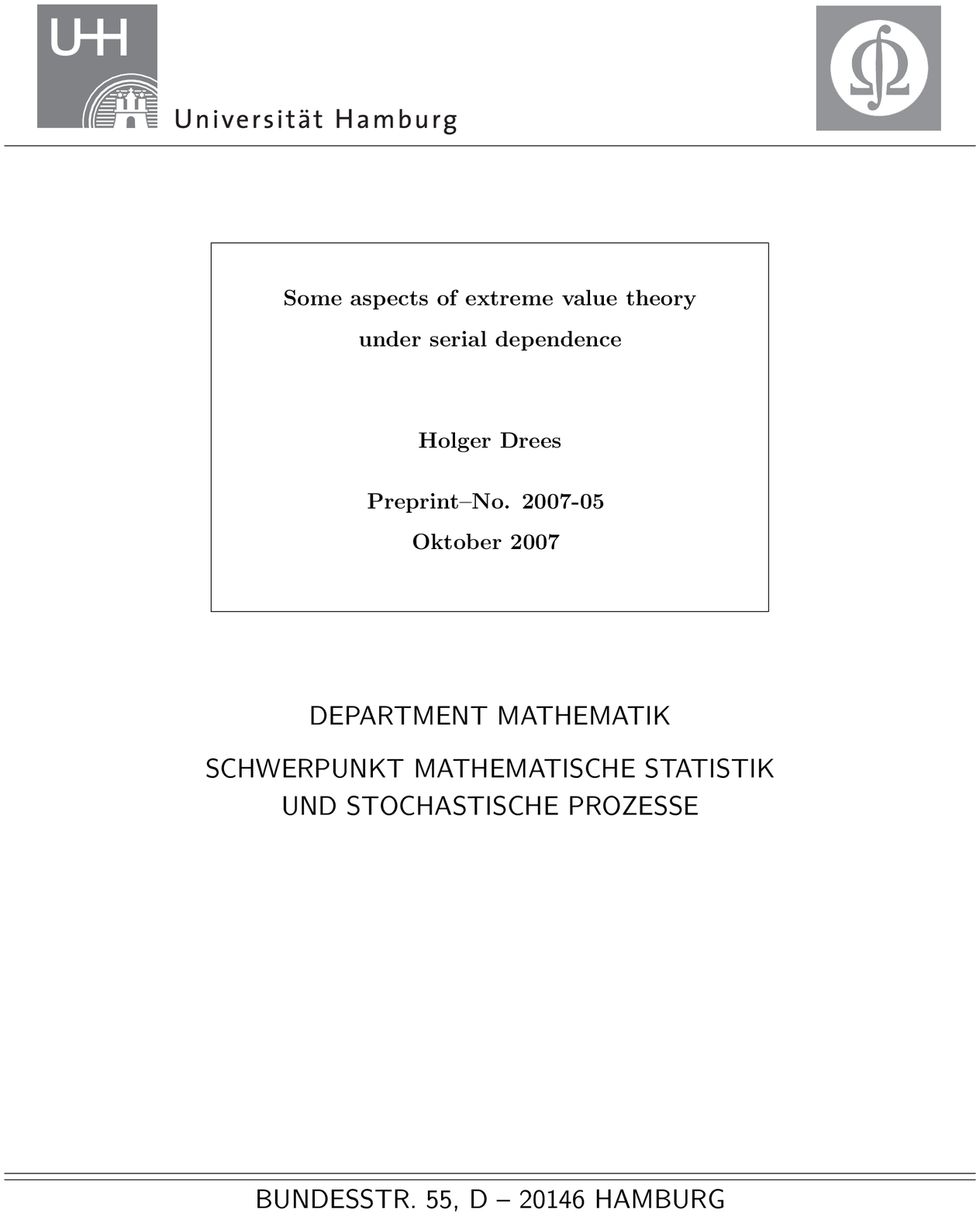}}

\vspace*{-4.5cm}

\setcounter{page}{0}

\maketitle

\centerline{\bf Abstract}

On the occasion of Laurens de Haan's 70th birthday, we discuss two aspects of
the statistical inference on the extreme value behavior of time series with a
particular emphasis on his important contributions. First, the performance of a
direct marginal tail analysis is compared with that of a model-based approach
using an analysis of residuals. Second, the importance of the extremal index as
a measure of the serial extremal dependence is discussed by the example of
solutions of a stochastic recurrence equation.
 \vspace{1ex}




\section{Introduction}

Since the publication of his Ph.D.-thesis, Laurens de Haan has been one of the main
driving forces behind the impressive development of extreme value statistics in the last
four decades. While he is best known for his seminal contributions to extreme value
theory for i.i.d.\ samples of univariate and multivariate observations and, in recent
years, for i.i.d.\ copies of continuous stochastic processes, he has also strongly
influenced the extreme value theory (and practice) for serially dependent data in two
ways: first by direct contributions, and second indirectly by promoting general
principles. In the present paper, both aspects of the impact of Laurens de Haan's work on
the development of extreme value statistics for time series are discussed.

Throughout his work, Laurens de Haan has always aimed at the greatest
(reasonable) generality of the models under consideration. For example, while
in many articles on univariate extreme value statistics it was assumed that
exact generalized Pareto random variables (r.v.s), respectively, generalized
extreme value r.v.s were observed, typically he only assumed that the
underlying distribution belongs to the domain of attraction of some extreme
value distribution. Under this much more general condition, he analyzed the
consequences of this deviation from the ideal situation. The second order
condition, de Haan and Stadtm{\"u}ller (1996) introduced and analyzed for that
purpose, is now {\em the} generally accepted standard condition in this field
(cf.\ \eqref{secordinno} for a simplified version). (It is worth mentioning
that essentially the same condition has been independently suggested by Pereira
(1994).) Similarly, also in de Haan's work on multivariate extreme value
statistics it is not assumed that the observations are drawn from an exact
extreme value distribution, but only that the underlying distribution belongs
to the domain of attraction of such a distribution. Moreover, he always
preferred weak smoothness conditions on the exponent measure pertaining to this
multivariate extreme value distribution to restrictive parametric submodels of
the natural infinite-dimensional extreme value model. With the improvement of
the resulting nonparametric methods, in the last couple of years this approach
has become more widely accepted as a reliable tool, that is more robust than
parametric approaches.

The extreme value estimators, that were suggested and analyzed for univariate i.i.d.\
samples by Laurens de Haan and many others, can also be used for the marginal tail
analysis of stationary time series, but often their performance deteriorates because of
the serial dependence between the observations. Therefore, in contrast to the
aforementioned general trend towards weak model assumptions, in the literature on extreme
value statistics for time series (and particularly linear time series) often an approach
is favored in which a parametric serial dependence structure is assumed and estimators
are considered which are based on a tail analysis of the (nearly independent) residuals
after the parametric model has been fitted. In the main Section 2, we will reassess some
of the results which seemingly show the underperformance of a direct extreme value
analysis, that only requires weak nonparametric assumptions on the dependence structure,
relative to the model-based approach.

Often one is interested not only in the marginal tail behavior but also in the
extremal dependence structure. The literature on the dependence analysis is
strongly dominated by the problem of estimating the extremal index, that
describes the influence of the serial dependence on the asymptotic behavior of
maxima of consecutive observations. It is somewhat surprising that, while in
the last two decades statistical methods which are based on exceedances (or
order statistics) instead of block maxima have become much more popular, this
shift of focus is not reflected in the statistical theory of the extremal
dependence structure. In Section 3, we will argue that the statistical
inference of the extremal dependence structure should be put on a broader basis
and exemplify this claim by the asymptotic behavior of naturally arising
statistics that were analyzed by Laurens de Haan and co-authors in a specific
time series model.

Obviously, the extreme value statistics of time series is a field of research
much too broad and diverse to be reviewed in a short article. For that reason,
we decided to focus on the two above topics, knowing that this selection is
largely a matter of taste. Important subfields which we will not discuss at all
are, for instance, the extreme value inference under additional structural
assumptions (e.g.\ for Markov chains) and the analysis of nonstationary or
multivariate time series, among many other topics. We will also not discuss
Laurens de Haan's contributions to the extreme value theory of continuous time
processes, since he usually assumes that i.i.d.\ copies of the whole process
are observed. Consequently, this theory is a natural extension of the theory
for multivariate observations rather than the theory for time series and will
thus be discussed in Michael Falk's contribution to this volume.

Throughout this article, we will assume that $X_t$, $t\in\Z$, is a stationary
time series with marginal distribution function (d.f.) $F_X$, that belongs to
the domain of attraction of some extreme value distribution.

\section{ Marginal tail analysis: In models we trust?}

In this section we assume that only the tail behavior of the stationary marginal d.f.\
$F_X$ is to be analyzed. To this end, estimators of tail parameters based on exceedances
over high thresholds can be used which were developed for i.i.d.\ observations, but the
serial dependence must be taken into account when the accuracy of these estimators is
assessed (e.g.\ to construct confidence intervals).

Roughly speaking, one can distinguish three different approaches:
\begin{enumerate}
  \item One tries to identify independent clusters of exceedances and constructs a new data set
by taking one observation (usually the cluster maximum) from each cluster. This way one
obtains an i.i.d.\ sample whose tail behavior can be analyzed using standard techniques
from classical extreme value statistics for i.i.d.\ data.
 \item In a nonparametric approach, one may also apply the classical tail estimators
(originally proposed for i.i.d.\ samples) directly to all exceedances observed in the
time series. However, if one wants to construct confidence intervals, then one needs
results on their asymptotic behavior that hold true under mild assumptions on the serial
dependence structural.
 \item Finally, in a semiparametric fashion, one can fit a parametric model of the serial
dependence to the data and then one can try to infer the tail behavior of the
time series from a suitable analysis of the residuals. This approach seems best
suited for heavy-tailed linear time series for which the relationship between
the tail of the stationary distribution and the tail of the distribution of the
innovations is particularly simple.
\end{enumerate}

The declustering approach (i) is most appropriate if the time series consists of clearly
separable, short clusters of extreme events, that preferably have a ``physical''
interpretation. Nice examples are data sets of wave heights and other quantities
describing sea conditions that were analyzed by Laurens de Haan and co-authors in several
publications. For example, starting with wave heights, wave periods and still water
levels that were observed every 3 hours at some point near the Dutch coast, de Haan and
de Ronde (1998) obtained nearly i.i.d.\ data by only considering the maximum of each
coordinate in each storm. See Dekkers and de Haan (1989), and de Haan (1990) for further
examples of that type.

Unfortunately, in many applications either it is difficult to identify {\em independent}
clusters of extremes, or the clusters are large so that it would be a great waste of
information to use but one observation in each cluster. For example, time series of
returns of some financial investment often exhibit long periods of high volatility during
which several dependent exceedances occur. Moreover, declustering schemes often depend on
certain subjective choices, and usually the influence of these choices on the accuracy of
the tail analysis is difficult to assess. (A first trial to overcome these problem was
made by Ledford and Tawn (2003).)

For these reasons, in the sequel we will focus on the nonparametric approach
(ii) and the model-based semiparametric approach (iii). In particular, we will
compare the accuracy and the robustness of resulting estimators of extreme
quantiles in the case of heavy-tailed linear time series.

\subsection{Direct extreme value analysis}

Among all tail estimators, the asymptotic behavior of the Hill estimator
 $$ \hat\gamma_{n,k} := \frac 1k \sum_{i=1}^n \log \frac{X_{n-i+1:n}}{X_{n-k:n}}
 $$
under serial dependence has been studied most thoroughly in literature; here
$X_{j:n}$ denotes the $j$th smallest order statistic of $X_1,\ldots,X_n$. One
of the first references is an unpublished manuscript by Rootz{\'e}n, Leadbetter and
de Haan (1990), in which the asymptotic normality of the Hill estimator is
established under quite weak conditions, including strong mixing of the time
series. At about the same time, Hsing (1991) independently proved the
asymptotic normality of the Hill estimator under comparable, but different
structural assumptions on the serial dependence. Since then, the limit
distribution of (variants of) the Hill estimator under serial dependence has
been examined in several papers; see, e.g., Resnick and St\u{a}ric\u{a} (1997)
and Novak (2002).

Of course, in most applications, the extreme value index is not the primary
object of interest, but for instance exceedance probabilities or extreme
quantiles are to be estimated. Consequently,  Rootz{\'e}n, Leadbetter and de Haan
(1990) also examined the asymptotic behavior of extreme quantile estimators.
Moreover, more general statistics of the type
 $$ \sum_{i=1}^n \phi_n\big( (X_i-u_n)^+\big) $$
(with suitable functions $\phi_n$) were considered, which are nowadays known as
tail array sums. The asymptotic theory of tail array sums was further developed
by Leadbetter and Rootz{\'e}n (1993) and Leadbetter (1995). In a final version,
this part of the manuscript was published in the article Rootz{\'e}n, Leadbetter
and de Haan (1998).

The general results about the asymptotic normality of tail array sums proved a powerful
tool. In particular, Rootz{\'e}n (1995), who established the weak convergence of the
empirical process
 $$ e_n(x) = \frac 1{\sqrt{n \bar F_X(u_n)}} \sum_{i=1}^n \Big( 1_{\{X_i>\sigma_n x+u_n\}}-
 \bar F_X(\sigma_n x+u_n)\Big)
 $$
 (with $\bar F_X=1-F_X$) towards a Gaussian  process under $\beta$-mixing (absolute regularity) of the time series,
  used this result to verify the
convergence of the finite dimensional marginal distributions. (In the improved version
Rootz{\'e}n (2006), a similar result is also established under the weaker assumption that the
time series is strongly mixing.) Using this convergence, Drees (2000) proved a weighted
approximation of the pertaining tail quantile process, from which one can easily conclude
the asymptotic normality of a general class of estimators of the extreme value index or
of estimators of extreme quantiles; cf.\ Drees (2000,2002,2003).

\subsection{ Model-based tail estimators}

 Here we focus on linear time series models, because
for this class the semiparametric  model-based approach seems particularly promising.
More precisely, we assume that the time series allows a representation as a moving
average of infinite order:
\begin{equation} \label{lints}
  X_t = \sum_{j=-\infty}^\infty \psi_j Z_{t-j}, \qquad t\in\Z.
\end{equation}
 Moreover, the i.i.d.\ innovations $Z_t$ are assumed to have balanced heavy
tails, i.e.\ their survival function $\bar F_Z$ satisfies
 \begin{equation} \label{heavybalance}
  \bar F_Z(x)=x^{-1/\gamma}L(x), \qquad \lim_{x\to\infty} \frac{\bar F_Z(x)}{\bar F_{|Z|}(x)}=p
 \end{equation}
for some $\gamma >0$, $p\in (0,1]$ and some slowly varying function $L$. Mikosch and
Samorodnitsky (2000) proved that
 \begin{equation} \label{lintailrel}
 \lim_{x\to\infty} \frac{\bar F_X(x)}{\bar F_{Z}(x)}= \frac 1p \sum_{j=-\infty}^\infty \Big( p\psi_j^{1/\gamma} 1_{\{\psi_j>0\}}
   + (1-p)|\psi_j|^{1/\gamma} 1_{\{\psi_j<0\}} \Big)
 \end{equation}
if $0<\sum_{j=-\infty}^\infty \psi_j^2<\infty$ for $\gamma<1/2$, and
$0<\sum_{j=-\infty}^\infty |\psi_j|^{1/\gamma-\eps}<\infty$ for some $\varepsilon
>0$ in the case $\gamma \geq 1/2$,
and if $E(Z_t)=0$ for $\gamma <1$. (Under stronger conditions, similar results were
already established by Davis and Resnick (1985) and Datta and McCormick (1998), among
others.)

 In particular, the time series
has the same extreme value index $\gamma$ as the innovations. Hence, if one has
estimated the coefficients $\psi_j$ and the time series model is invertible,
then it suggest itself to estimate $\gamma$ by applying the Hill estimator (or
some other estimator of the extreme value index) to the resulting residuals
$\hat Z_t$, which are approximately i.i.d. Moreover, if $p$ (or some estimate
of it) is known, one may even calculate estimates of exceedance probabilities
$\bar F_X(x)$ over high thresholds $x$ or estimators of extreme quantiles
$F_X^{-1}(1-t)$ (for small $t>0$) from estimators of the corresponding
quantities of the d.f. $F_Z$ of the innovations, which in turn can be obtained
from a tail analysis of the residuals. This program has been worked out for the
first time by Resnick and St\u{a}ric\u{a} (1997) for the Hill estimator and
autoregressive time series $AR(m)$ of order $m<\infty$:
 $$X_t=\sum_{i=1}^m \varphi_i X_{t-i}+ Z_t, \quad t\in\Z.$$
 Let $\hat \varphi_i$, $1\le i\le m$, be estimators
of the coefficients such that
 $d_n(\hat \varphi_i-\varphi_i)_{1\le i\le m}$ converges to some nondegenerate
distribution; here
 $d_n \rightarrow \infty$ determines the rate of convergence of the
estimators $\hat \varphi_i$. Define the residuals
 $$\hat Z_t:=X_t-\sum_{i=1}^m \hat\varphi_i X_{t-i}, \qquad m+1\leq
t\leq n.$$
 Resnick and St\u{a}ric\u{a} (1997) proved that the Hill estimator based on the
absolute residuals
$$\hat\gamma_{n,|Z|}:= \frac 1{k} \sum_{i=1}^{k} \log \frac{|\hat Z|_{n-m-i+1:n-m}}{|\hat Z|_{n-m-k:n-m}}$$
 (with $|\hat Z|_{j:n-m}$
denoting the $j$th smallest order statistic of $|Z_t|$, $m-1\leq t\leq n$) is
asymptotically normal
 $$ \sqrt{k} \big( \hat\gamma_{n,|Z|}-\gamma\big) \; \longrightarrow \mathcal{N}_{(0,\gamma^2)} $$
 weakly, provided that the d.f.\ $F_{|Z|}$ of the absolute
innovations satisfies the second order condition
 \begin{equation} \label{secordinno}
    \lim_{t\to\infty} \frac{\displaystyle \frac{\bar F_{|Z|}(tx)}{\bar F_{|Z|}(t)} x^{1/\gamma} -1}{\tilde
   A(t)} = \frac{x^{\tilde\rho}-1}{\tilde\rho}
 \end{equation}
for some $\tilde\rho\le 0$, and the number $k$ of order statistics used for estimation
tends to $\infty$ sufficiently slowly such that
 \begin{equation} \label{kcond}
  \lim_{n\to\infty} \sqrt{k} \tilde A \big(F^{-1}_{|Z|} (1 - k/n)
\big) =  0
 \end{equation}
 and
 \begin{equation} \label{extrakcond}
   \lim_{n\to\infty} \frac{\displaystyle\sqrt{k} F^{-1}_{|Z|} (1 - \sqrt{k}/n)}{\displaystyle d_n F^{-1}_{|Z|} (1 - k/n)}
 = 0.
 \end{equation}
(Ling and Peng (2004), who established a similar result for ARMA time series, showed that
condition \eqref{extrakcond} is not needed if $d_n$ equals the best attainable rate.)

In contrast, the Hill estimator applied to the absolute observations $|X_t|$ directly is
asymptotically normal with asymptotic variance
 \begin{equation} \label{asvardir}
   \gamma^2 \bigg( 1+2\frac{\sum_{j=1}^\infty \sum_{i=0}^\infty \min\big(
   |\psi_j|^{1/\gamma},|\psi_{i+j}|^{1/\gamma}\big)}{\sum_{i=0}^\infty
   |\psi_j|^{1/\gamma}}\bigg)
 \end{equation}
 with $\psi _i$ denoting the coefficients
of the $MA(\infty )$-representation of the time series, i.e.
 $ X_t=\sum_{i=0}^\infty \psi_i Z_{t-i}$, provided that the d.f.\ $F_{|X|}$ of the absolute
observations and the sequence of numbers of order statistics used for
estimation satisfy the analogs to the conditions \eqref{secordinno} and
\eqref{kcond}. Note that the asymptotic variance \eqref{asvardir} of the Hill
estimator directly applied to the absolute values of the observations is
strictly larger than the asymptotic variance $\gamma ^2$ of the Hill estimator
based on the absolute residuals. For example, if one considers an $AR(1)$ time
series with coefficient $\varphi_1\in (-1,1)$, then $\psi _i=\varphi_1^i$,
$i\in\N_0$, and the asymptotic variance \eqref{asvardir} equals $\gamma
^2(1+|\varphi_1|^{1/\gamma})/(1-|\varphi_1|^{1/\gamma})$. Therefore, Resnick
and St\u{a}ric\u{a} (1997) claimed that ``the procedure of applying the Hill
estimator directly to an autoregressive process is less efficient than the
procedure of first estimating autoregressive coefficients and then estimating
$\alpha$'' ($=1/\gamma$) ``using estimated residuals''. This conclusion,
however, is justified in general {\em only if both Hill estimators use the same
number of order statistics}. Since the optimal numbers of order statistics used
by the Hill estimators are essentially determined by the functions $\tilde A$
occurring in the second order condition \eqref{secordinno} for $F_{|Z|} $ (in
the case of the residual-based estimator) and $A$ in the analog condition for
the absolute time series (for the directly applied Hill estimator), it is
 a priori unclear which of the estimators has the smaller variance
when they both use an appropriate number of order statistics. Indeed, if the
second order parameter$\tilde\rho$ is smaller in absolute value than the
analogous parameter $\rho$ from the second order condition for $F_{|X|}$, then
the best rate of convergence that can be achieved by the residual-based Hill
estimator is of lower order than the optimal rate of the directly applied Hill
estimator, i.e.\ the former estimator has asymptotic efficiency 0 with respect
to the latter estimator. Conversely, if $|\tilde\rho| >|\rho|$ then the
directly applied Hill estimator is asymptotically inefficient w.r.t.\ the
model-based estimator, if both use the optimal number of order statistics.

 For general linear time series, it is not known how the second
order behavior of $F_{|X|} $ is related to the second order behavior of
$F_{|Z|}$. However, for first order moving averages the relationship has been
discussed by Geluk, Peng and de Vries (2000) and Geluk and Peng (2000), and the
same technique can be used for finite order moving averages. For general linear
dependence structures but a rather particular class of distributions of
innovations, the relationship can be deduced from results by Barbe and
McCormick (2004). More precisely, assume that
$$ \bar F_Z(x) = x^{-1/\gamma}(c+d x^{-1}+o(x^{-1})), \quad F_Z(-x) = x^{-1/\gamma}(\tilde c+\tilde d
x^{-1}+o(x^{-1}))
$$
as $x\to\infty$. Then, by Section 2.1 of Barbe and McCormick (2004), the tail
of the linear time series \eqref{lints} behaves asymptotically as
\begin{eqnarray*}
 F_X(x) & = & x^{-1/\gamma}(d_\psi+ D_\psi x^{-1}+o(x^{-1})) \quad
 \text{with}\\
 d_\psi & = & \sum_{j=-\infty}^{\infty} \Big( c\psi_j^{1/\gamma} 1_{\{\psi_j>0\}}
 + \tilde c|\psi_j|^{1/\gamma} 1_{\{\psi_j<0\}}\Big), \\
 D_\psi & = & \sum_{j=-\infty}^\infty \Big( cd\psi_j^{1/\gamma+1} 1_{\{\psi_j>0\}}
 + \tilde c\tilde d|\psi_j|^{1/\gamma+1} 1_{\{\psi_j<0\}}\Big).
\end{eqnarray*}
Hence, for this type of innovations,  both functions $A(t)$ and $\tilde A(t)$
are multiples of $t^{-1}$, but the constant factors differ from each other.
Note that the above expansion of $\bar F_Z$ is equivalent to
$F_Z^{-1}(1-t)=c^\gamma t^{-\gamma} +\gamma d+o(1)$, i.e., up to terms of
smaller order the tails behave like those of a shifted Pareto distribution.
This shows that indeed the result by Barbe and McCormick describes the
relationship between the second order behavior of the tails of $F_Z$ and $F_X$
(or of $F_{|Z|}$ and $F_{|X|}$) only for a rather limited family of
distributions of innovations.

To sum up, in general, the result by Resnick and St\u{a}ric\u{a} (1997) does
not allow to compare the asymptotic performance of the direct nonparametric
approach and the model-based estimator if both use an appropriate number of
order statistics.

Having said that, it is nevertheless plausible to expect that the residual-based
estimator has a smaller variance if the extra factor by which $\gamma ^2$ is multiplied
in the variance formula \eqref{asvardir} is much larger than 1 (e.g.\ if the absolute
coefficient of an $AR(1)$ time series is close to 1). However, even in that case, the
model based approach has serious drawbacks:
\begin{itemize}
  \item As mentioned before, usually one is not mainly interested in
the extreme value index but e.g.\ in exceedance probabilities $\bar F_X(x)$ (or extreme
quantiles). To estimate such quantities, one uses the relationship \eqref{lintailrel}
which may introduce a non-negligible additional error term if for the given threshold $x$
the ratio $\bar F_X(x)/\bar F_Z(x)$ is poorly approximated by the right hand side of
\eqref{lintailrel}.
 \item Of course, the model-based approach makes sense only if the
model assumptions are (approximately) fulfilled. Although this remark is almost
trivial, it is nevertheless crucial to be aware of the fact that even moderate
deviations from the linear relationship between $X_t$ and $X_{t-1},\ldots,
X_{t-m}$, which can hardly be detected by statistical tests, may completely
wreck up the residual-based estimator, as we will see in the next subsection.
\end{itemize}

\subsection{ Comparison of model-based and direct extreme quantile estimators: a simulation
study}

 Here we assume that $X_t$, $t\in \mathbb{Z}$, is a stationary $AR(1)$ time series,
i.e. $X_t=\varphi_1 X_{t-1}+Z_t$ with i.i.d.\ innovations $Z_t$ satisfying
\eqref{heavybalance}, and that an extreme quantile $F^{-1}_X(1-t)$ ($t>0$
small) is to be estimated. In this case relation \eqref{lintailrel} reads as
 \begin{equation} \label{lintailrelAR1}
    \lim_{x\to\infty} \frac{\bar F_X(x)}{\bar F_Z(x)} = \left\{
     \begin{array}{lcl}
        1/(1-\varphi_1^{1/\gamma}), & & \varphi_1\in [0,1),\\
        \big( 1+|\varphi_1|^{1/\gamma}(1-p)/p\big)/(1-|\varphi_1|^{2/\gamma}), & &
        \varphi_1\in (-1,0).
     \end{array}
     \right.
 \end{equation}
For simplicity (and in favor of the model-based approach) we assume that $p$ is
known to be equal to $1/2$, so that \eqref{lintailrelAR1} simplifies to
 $$\lim_{x\to\infty} \frac{\bar F_X(x)}{\bar F_Z(x)} = 1/(1-|\varphi_1|^{1/\gamma}), $$
  which, by the regular
variation of $\bar F_Z$, is equivalent to the following relationship between the
corresponding quantile functions:
 $$ \lim_{t\downarrow 0} \frac{F_X^{-1}(1-t)}{F_Z^{-1}\big(1-(1-|\varphi_1|^{1/\gamma})t\big)}=1. $$
 Hence, in the model-based
approach one may estimate $F^{-1}_X (1-t)$ as follows:
 \begin{itemize}
  \item Estimate $\varphi_1$, e.g., by the sample auto-correlation $\hat\varphi_1$ at lag 1.
  \item Approximate $F^{-1}_X (1-t)$ by $F^{-1}_Z \big(1-(1-|\hat\varphi_1|^{1/\gamma})t\big)$
  and estimate the latter by
the Weissman estimator
 $$\hat Z_{n-k:n-1} \bigg( \frac{n(1-|\hat\varphi_1|^{1/\hat\gamma_{n,Z}})t}k\bigg)^{\hat\gamma_{n,Z}}
 $$
  where $\hat Z_{j:n-1}$ is
the $j$th smallest order statistic of the residuals $\hat Z_t=X_t-\hat\varphi_1 X_{t-1}$,
$2\leq t\leq n$, and
 $$ \hat \gamma_{n,Z}:= \frac 1k \sum_{i=1}^k \log \frac{\hat Z_{n-i:n-1}}{\hat
 Z_{n-k-1:n-1}}
 $$
  is the corresponding Hill estimator. (The Weissman
estimator can be motivated either by a Generalized Pareto approximation of $F_Z$ or the
regular variation of $F^{-1}_Z$ which implies $F^{-1}_Z (1-u) \approx
F^{-1}_Z(1-k/n)(nu/k)^{-\gamma}$ for sufficiently small $u$ and $k/n$.)
 \end{itemize}

In a small simulation study we consider time series of length $n=2000$
 with $\varphi_1=0.8$ and two different symmetric distributions of
the innovations $Z_t$:
 \renewcommand{\labelenumi}{\alph{enumi})}
\begin{enumerate}
 \item a double-sided (unshifted) Pareto distribution, i.e. $\bar F_Z(x)=F_Z(-x)=0.5 x^{-1/\gamma}$ for $x\ge 1$ with
 $\gamma=1/2$;
 \item a double-sided shifted Pareto distribution, i.e. $\bar F_Z(x)=F_Z(-x)=0.5(x+1)^{-1/\gamma}$ for $x\ge 0$  with
 $\gamma=0.3$.
\end{enumerate}
Clearly, Model a) is particularly favorable for the model-based approach, since
the Hill estimator based on the innovations and (according to aforementioned
results) also the Hill estimator based on the residuals is asymptotically
unbiased for all intermediate sequences $k=k_n$. In contrast, one might expect
a significant bias of the model-based quantile estimator in Model b) if one
uses too large a number of order statistics, as the Hill estimator is sensitive
to a shift of the data. In the first Model, it is {\em a priori} not clear how
large the bias of the direct nonparametric quantile estimator will be, since
the second order behavior of $F_X$ is not known.  For Model b), however, the
aforementioned result  by Barbe and McCormick (2004) is applicable. A careful
inspection of the proofs given by Resnick and St\u{a}ric\u{a} (1997) and Drees
(2000) and lengthy, but straightforward calculations show that the ratio of
minimal asymptotic root mean squared errors of the model-based Hill estimator
and the directly applied Hill estimator equals
$$ \bigg(
\frac{(1-|\varphi_1|^{1/\gamma+1})^2}{(1-|\varphi_1|^{1/\gamma})^2(1+|\varphi_1|^{1/\gamma})^{2\gamma}}\bigg)^{1/(2\gamma+1)},
$$
which is approximately equal to 1.03 for $\varphi_1=0.8$ and $\gamma=0.3$.
Hence, one may expect that the directly applied Hill estimator performs
slightly better than the model-based estimator and this superiority may also
carry over to the corresponding quantile estimators.

 We assume that the quantile
$F_X^{-1}(0.999)$ is to be estimate. In both models, this quantile is
approximated by the average of the corresponding empirical quantiles of 200
simulated time series of length $9\cdot 10^6$ (such that $F_X^{-1}(0.999)$ lies
well within the simulated data sets), which yields $F_X^{-1}(0.999)\approx
37.94$ in Model a) and $F_X^{-1}(0.999)\approx 7.312$ in Model b) (with a
relative approximation error of less than $0.002$ with probability of at least
0.99).

Figure 1 displays the simulated root mean squared error (RMSE) and the
$L_1$-error of the direct quantile estimator (solid, resp.\ dotted line) and of
the model-based estimator (dashed, resp.\ dash-dotted line) versus the number
$k$ of order statistics used for estimation. Obviously, the model-based
approach outperforms the direct estimator in Model a), in that it has a much
smaller RMSE and $L_1$-error if $k$ is chosen optimally. Moreover, its
performance is less sensitive to an inappropriate choice of $k$: it performs
reasonably well for all values of $k$ between 150 and 750 (i.e., as expected
one might use a very large proportion of all positive residuals), while the
performance of the direct quantile estimator quickly deteriorates when $k$ is
smaller than 200 or larger than 300. Conversely, as expected, in Model b) the
direct estimator performs somewhat better than the model-based estimator, i.e.,
its minimal RMSE and $L_1$-error is smaller than the corresponding errors of
the model-based estimator, and the performance of the direct estimator is less
sensitive to the choice of $k$. However, both effects are much less pronounced
than in Model a). This can also be seen from Table 1 which summarizes the
minimal errors with the corresponding optimal values of $k$ and also the
simulated bias and the standard errors (i.e.\ simulated standard deviations)
for the choice of $k$ which leads to the minimal RMSE: while in Model a) the
RMSE of the direct estimator is about 2.5 times larger than the RMSE of the
model-based estimator, the ratio between the minimal RMSE's in Model b) is just
about 1.2.

\begin{figure}


\includegraphics[height=50mm]{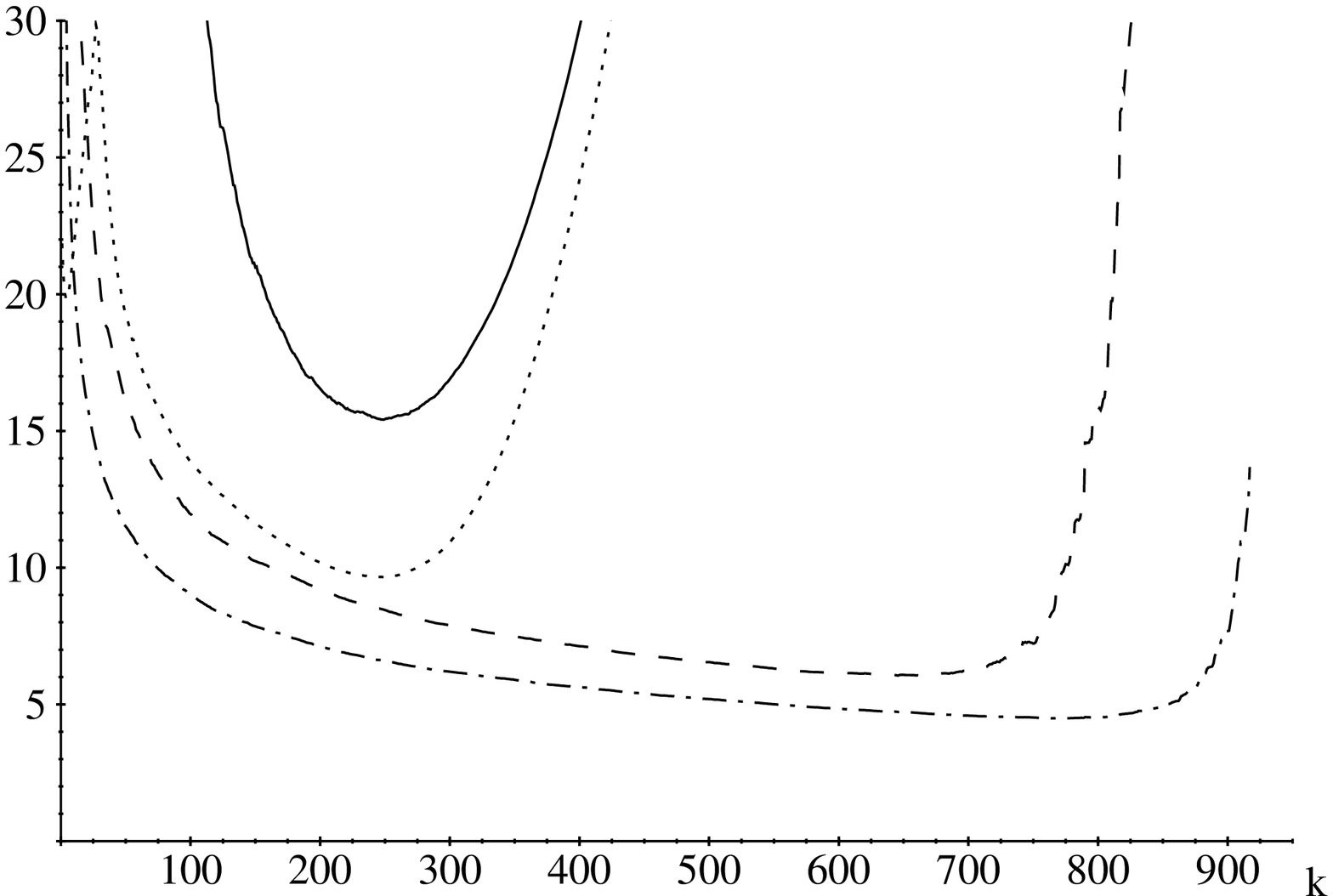} \hfill
\includegraphics[height=50mm]{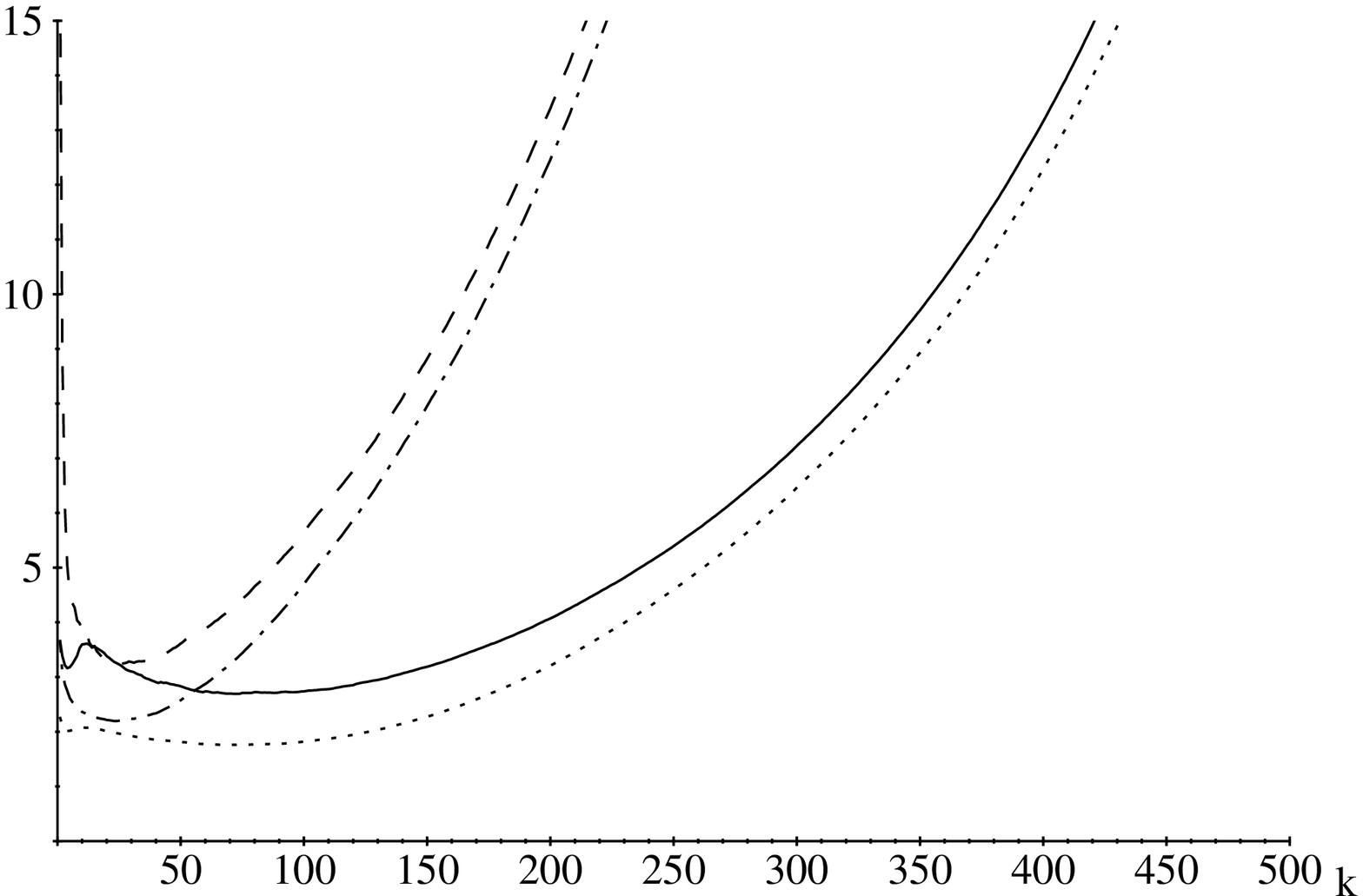}

\caption{Simulated RMSE and $L_1$-error of the quantile estimator directly applied to the
data (solid resp.\ dotted line) and of the quantile estimator based on the analysis of
residuals (dashed resp.\ dash-dotted line) vs.\ number $k$ of order statistics for a
linear $AR(1)$ time series; left plot: innovations according to Model a), right plot:
Model b) }
\end{figure}

 \vspace{2ex}

\begin{table}
\begin{tabular}{l|cc|cc}
 & \multicolumn{2}{c}{Model a)} &  \multicolumn{2}{c}{Model b)}\\
 &  RMSE  &$L_1$-error  &  RMSE  &$L_1$-error  \\
 & bias / s.e. & & bias / s.e.  \\ \hline
direct estimation & 15.4 \;(k=249) & 9.7 \;(k=249) & 2.7\; (k=73)    &   1.8 \;(k=70) \\
 & 3.3 / 15.0 & & 0.7 / 2.6\\[1ex]
model based estimation &  6.1 \;(k=662) & 4.5 \;(k=765) & 3.2 \; (k=25)  &     2.2 (k=23)\\
 & -0.6 / 6.0 & & 1.1 / 3.0
\end{tabular}
\caption{Minimal errors, bias and standard errors of the quantile estimators in the
(unperturbed linear) $AR(1)$ time series models}
\end{table}

From these results, one gets the impression that, although the model-based
approach does not yield more accurate estimators for all distributions of
innovations, its overall performance is  at least as good as the performance of
the direct nonparametric approach. However, as we will see next, this
interpretation is premature (and indeed quite dangerous) as the model-based
estimator can be very sensitive to relatively small deviations to the model.

As an example, we consider a nonlinear $AR(1)$ time series, namely a stationary solution
to the equation
\begin{equation}  \label{AR1pert}
 X_t=\varphi_1 X_{t-1}+ \delta \mathrm{sgn}(X_{t-1}) \log\big(\max(|X_{t-1}|,1)\big) +
Z_t \qquad \text{with } \varphi=0.8,\, \delta=0.6.
 \end{equation}
 Here the linear relationship between $X_t$ and
$X_{t-1}$ is perturbed by an extra logarithmic term. Of course, one cannot expect that
the relationship \eqref{lintailrelAR1} holds in this nonlinear model. Hence, most likely,
the model-based estimator will show an increased bias.

From Figure 2, which shows a typical scatterplot of $(X_{t-1},X_t)$ for a time
series according to model \eqref{AR1pert} with shifted double-sided Pareto
innovations from Model b), it is apparent that, with the naked eye, such a time
series can hardly be distinguished from a classical linear $AR(1)$ time series
(with an increased autoregressive coefficient). Moreover, if one fits a linear
$AR(1)$ model to such a time series of length $n=2000$, then the turning point
test and the difference-sign test with nominal size 0.05 (see Brockwell and
Davis (2002), p.\ 36 f.) detect dependence in the residuals with probability
less than 0.06, i.e.\ these tests are not capable of detecting the model
deviation. Also the Portmanteau test with nominal size 0.05 applied to the
residuals has a maximal power of about 0.13, i.e.\ the rejection rate is less
than 0.1 higher than the false alarm rate if the data comes from the
corresponding $AR(1)$ model. (Note that the Portmanteau test should not be
applied if the variance of the innovations is infinite; hence, strictly
speaking, it is not suitable for Model a).) To sum up, it is almost impossible
to distinguish a time series from model \eqref{AR1pert} from a linear $AR(1)$
time series.

\begin{figure}


\centerline{\includegraphics[height=50mm]{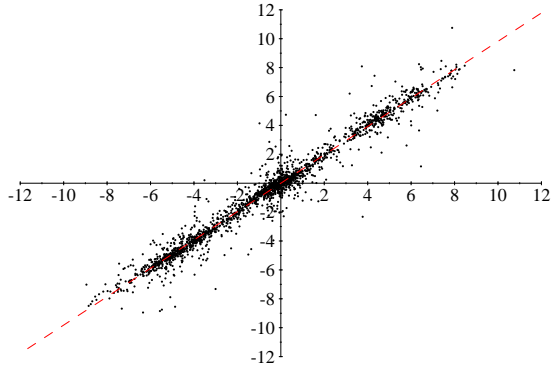}}

\caption{Scatterplot $(X_{t-1},X_t)$ of a simulated nonlinear $AR(1)$ time series
\eqref{AR1pert} of length $n=2000$ with innovations according to Model b); the dashed
line indicates the fitted linear relationship with estimated autoregressive coefficient
0.982}
\end{figure}

Figure 3 shows the RMSE and $L_1$-errors of the quantile estimator obtained from the
direct approach and the model-based approach when the time series are simulated according
to \eqref{AR1pert}, but erroneously a linear $AR(1)$ model is assumed. Table 2 gives the
minimal errors in this case (analogously to Table 1). For both d.f.s of the innovations,
the errors of the model-based quantile estimators are much larger in the nonlinear
$AR(1)$ model than for the classical linear $AR(1)$ time series. To a large extent, the
deterioration of the performance is caused by the large bias, but also the variance is
much larger now even if the decrease of the optimal number of order statistics is taken
into account. In sharp contrast, the direct quantile estimator, which does not rely on a
specific time series model, is more precise for these nonlinear $AR(1)$ time series than
for the linear ones. Consequently, if the innovations are drawn from Model a), then the
minimal RMSE of the model-based quantile estimator is about 25\% larger than the minimal
RMSE of the direct quantile estimator, while for innovations according to Model b) the
RMSE of the model-based estimator is more than 7 times larger than the RMSE of the
nonparametric estimator!

\begin{figure}


\includegraphics[height=50mm]{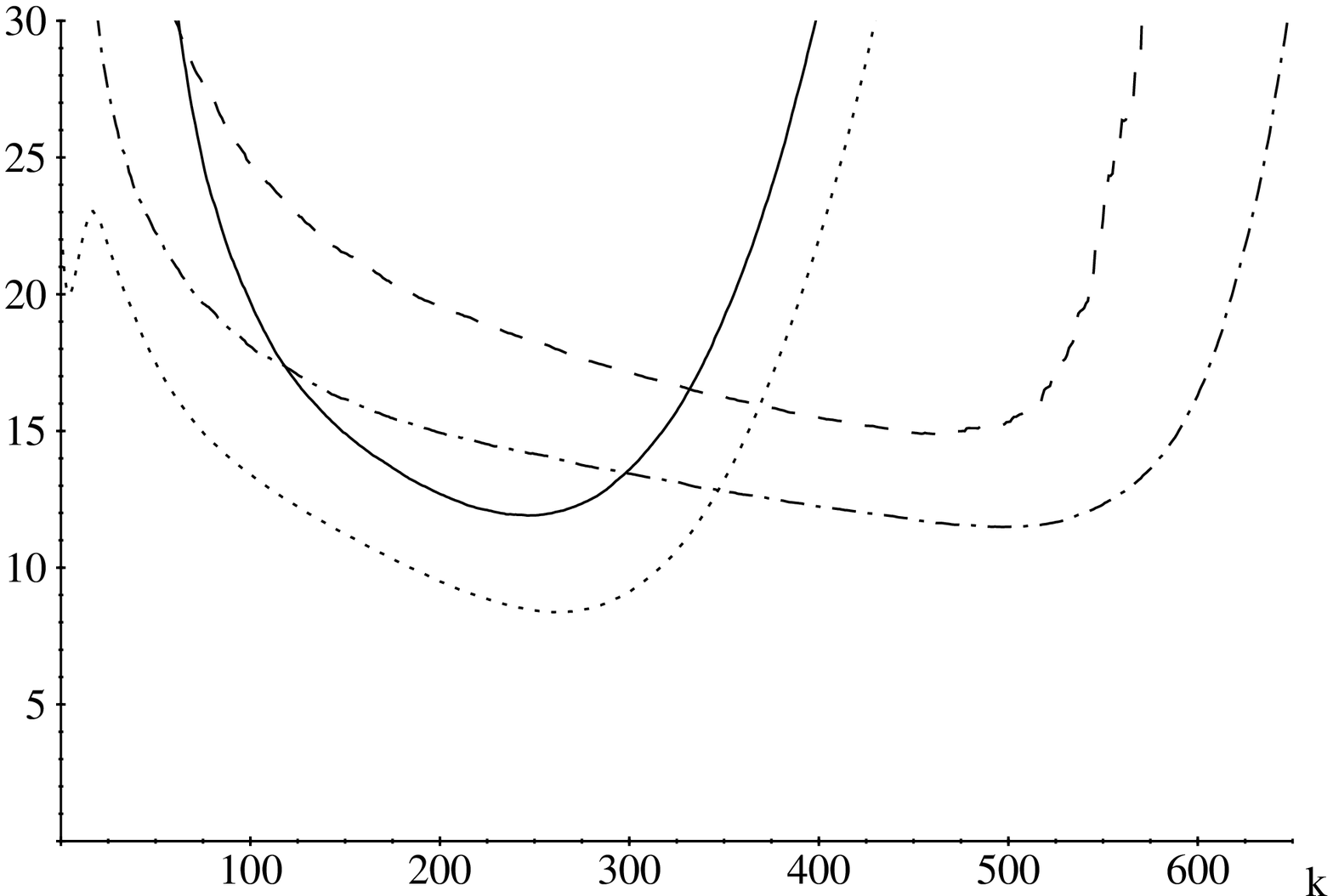} \hfill
\includegraphics[height=50mm]{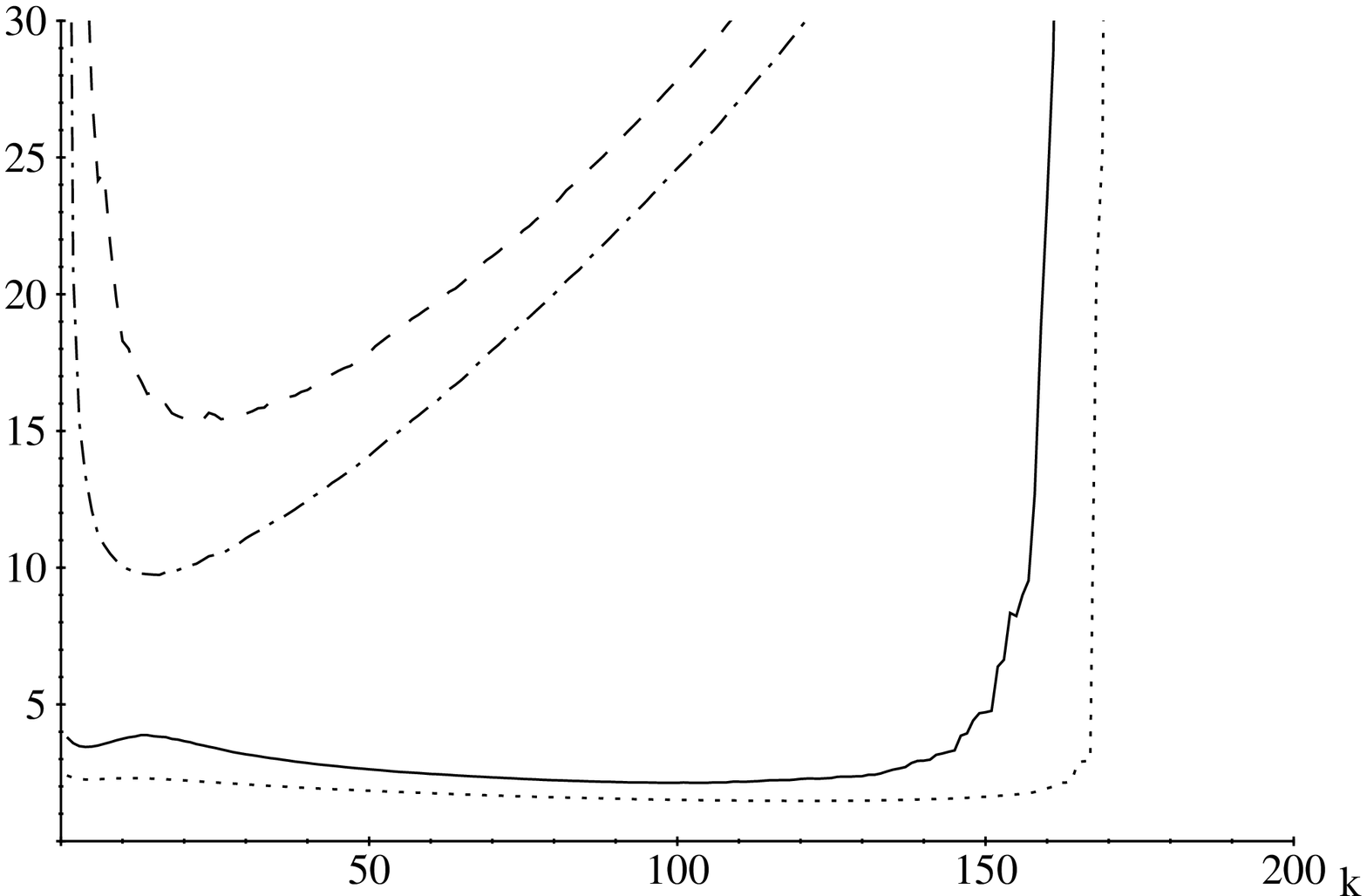}

\caption{Simulated RMSE and $L_1$-error of the quantile estimator directly applied to the
data (solid resp.\ dotted line) and of the quantile estimator based on the analysis of
residuals (dashed resp.\ dash-dotted line) vs.\ number $k$ of order statistics for the
nonlinear $AR(1)$ time series \eqref{AR1pert}; left plot: innovations according to Model
a), right plot: Model b) }
\end{figure}

 \vspace{2ex}

\begin{table}
\begin{tabular}{l|cc|cc}
 & \multicolumn{2}{c}{Model a)} &  \multicolumn{2}{c}{Model b)}\\
 &  RMSE  &$L_1$-error  &  RMSE  &$L_1$-error  \\
 & bias / s.e. & & bias / s.e.  \\ \hline
direct estimation &  11.9 \; (k=247)    &      8.4\; (k=260)&     2.1 \; (k=99) & 1.5 \;
(k=122) \\
 & -0.0 / 11.9 & & 0.0 / 2.1\\[1ex]
model based estimation &  14.9 \; (k=462)    &     11.4 \; (k=498) &  15.4\; (k=22) &        9.7\; (k=16)\\
 & 9.2 / 11.7 & & 9.5 / 12.1
\end{tabular}
\caption{Minimal errors, bias and standard errors of the quantile estimators in the
nonlinear $AR(1)$ time series \eqref{AR1pert}}
\end{table}

Figure 4 demonstrates that the very poor performance of the quantile estimator
which is based on the residual analysis is not due a few particularly wrong
estimates but that indeed the estimator yields rather poor results with a high
probability. In this plot, for innovations according to Model b), kernel
estimates of the density of the direct (solid line) and the model-based
quantile estimators (dashed line) are displayed for optimal values of $k$
(i.e., $k=99$ and $k=22$, respectively). While the mode of both densities is
close to the true value (indicated by the vertical dotted line), the density of
the model-based estimator is strongly skewed to the right and has a large
spread. In contrast, the distribution of the nonparametric estimator is more
symmetric and much more concentrated around the true value.

\begin{figure}

\vspace{-0.5cm}

\centerline{\includegraphics[height=50mm]{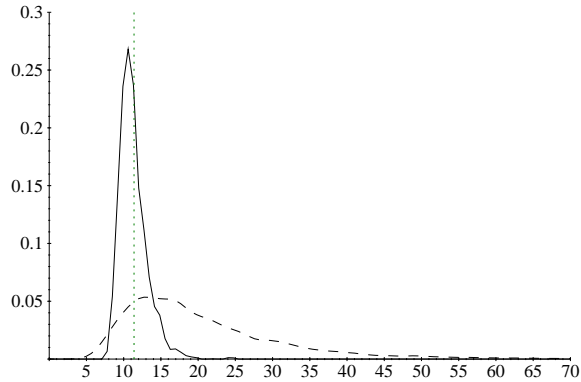}}

\caption{Estimated density of the direct quantile estimator (solid line) and the
model-based estimator (dashed line) for time series according to \eqref{AR1pert} with
innovations according to Model b); the true quantile is indicated by the vertical dotted
line}
\end{figure}

To sum up, the example shows that the model-based approach to the estimation of
extreme quantiles can give completely misleading estimates even if the
deviation from the assumed linear time series model is moderate in the sense
that it is very difficult to detect by means of statistical tests. Therefore,
it seems advisable to use estimators  only with utmost care which are based on
an extreme value analysis of residuals obtained under parametric model
assumptions about the dependence structure. In particular, it is not justified
to consider them generally superior to the directly applied extreme value
estimators.

\section[Analysis of the extremal dependence structure]{Analysis of
the extremal dependence structure: Is there a world behind the extremal index?}

So far we have only considered estimators of the marginal tail behavior. For many
applications, also the dependence between consecutive extreme values of the time series
is of interest. For example, if $X_t$ denotes the negative return (loss) of some
financial investment, it is important to assess the risk that all (or some of) the losses
$X_t,X_{t+1},\ldots,X_{t+m-1}$ in $m$ consecutive periods (or perhaps the total loss
$\sum_{i=0}^{m-1} X_{t+i}$) are large.

In the analysis of the extreme value behavior of maxima $M_n:=\max_{1\le t\le n} X_t$ of
$n$ consecutive observations the so-called extremal index plays a crucial role. Let
$\tilde X_t$, $1\le t\le n$, denote an associated sequence of i.i.d.\ random variables
with d.f.\ $F_X$. Assume that, for some normalizing constants $a_n>0$ and $b_n\in\R$,
$$ \frac{\max_{1\le t\le n} \tilde X_t-b_n}{a_n} \;\longrightarrow\; G $$
weakly for some nondegenerate d.f.\ $G$. Leadbetter (1983) proved that then
$$ \frac{M_n-b_n}{a_n} \;\longrightarrow\; G^\theta $$
for some $\theta\in [0,1]$ provided Leadbetter's condition $D(u_n)$
\begin{eqnarray*}
  \lefteqn{\hspace*{-1cm}\forall\, m,n>0 \;
  \exists \,\alpha_{n,m} \; \forall\, k,l>0, \, 1\le i_1<i_2<\ldots<i_k<i_k+m\le j_1<j_2<\ldots<j_l\le n:}\\
   & &\bigg| P\Big\{ \max_{1\le r\le k} X_{i_r} \le u_n,\; \max_{1\le s\le l} X_{i_s} \le
  u_n\Big\} - \\
  & & \hspace{3cm} - P\Big\{ \max_{1\le r\le k} X_{i_r} \le u_n\Big\}\cdot P\Big\{ \max_{1\le s\le l} X_{i_s} \le
  u_n\Big\}\bigg| \le \alpha_{n,m}\\[0.5ex]
  & \text{and} & \exists\, m_n=o(n): \quad \lim_{n\to\infty} \alpha_{n,m_n}=0
\end{eqnarray*}
holds for all $u_n=a_nx+b_n$ with $x>G^{-1}(0)$, and the d.f.\ $P\{(M_n-b_n)/a_n\le x\}$
of the standardized maximum converges for some $x>G^{-1}(0)$. Hence, if the extremal
index $\theta$ is strictly positive, the maximum converges to some nondegenerate limit
distribution that is of the same type as the limit distribution in the case of
independence.

Moreover, Hsing, H{\"u}sler and Leadbetter (1988) proved that under weak additional
assumptions (including the slightly stronger mixing condition $\Delta(u_n)$)
the point process $\sum_{t=1}^n \eps_{t/n} 1_{(a_nx+b_n,\infty)} (X_t)$ of
standardized time points at which exceedances occur converges to a compound
Poisson process. Then, typically, the extremal index $\theta$ equals the
reciprocal value of its mean cluster size (although in general one only knows
that $\theta$ is a lower bound for this value).

Since the asymptotic behavior of maxima of consecutive observations is completely
determined by the extremal index and the tail behavior of $F_X$, the literature on the
statistical analysis of the extremal dependence structure focusses on the estimation of
$\theta$ and, to a lesser extent, the estimation of the cluster size distribution; see
Hsing (1991, 1993), Smith and Weissman (1994), Weissman and Novak (1998),
Ancona-Navarrete and Tawn (2000) and Ferro and Segers (2003), among others.

However, as the aforementioned example shows, in financial applications often
other statistics of extreme values (than maxima) are of main interest, and the
same holds true in other fields of applications where exceedances over high
thresholds rather than block maxima are considered. We will demonstrate by a
particular time series model that the extremal index and the cluster size
distributions are often not sufficient to determine the distribution of
statistics which arise in a natural way.

In the remainder of this section, we consider stationary solutions of stochastic
recurrence equations of the type
\begin{equation} \label{receq}
 X_t=A_t X_{t-1}+B_t, \qquad t\in\Z,
\end{equation}
with $(A_t,B_t)$ denoting  i.i.d.\ random vectors with values in
$(0,\infty)^2$. For instance, a squared $ARCH(1)$ time series satisfies this
relationship; further applications of this model were described by Vervaat
(1979). Kesten (1973) proved that such a stationary solution exists if $A_1$
does not have a lattice distribution, the distribution of $B_1/(1-A_1)$ is not
degenerate, and if there exists $\kappa>0$ such that $E A_1^\kappa=1$,
$E(A_1^\kappa\max(\log A_1,0))<\infty$ and $EB_1^\kappa\in (0,\infty)$.
Moreover, then $\bar F_X(x)\sim cx^{-\kappa}$ for some $c>0$ so that the
standardized maxima of an accompanying i.i.d.\ sequence converge to a Fr{\'e}chet
distribution with extreme value index $\gamma=1/\kappa$.

De Haan et al.\ (1989) calculated the extremal index and the cluster size
distribution of such a time series. Let $W_j=\prod_{i=1}^j A_i^\kappa$ (with
the convention $W_0=1$). Note that $\log W_j, j\ge 1,$ is a random walk with
negative drift. Hence the sequence $W_j, j\ge 1,$ tends to 0 and it is almost
surely bounded. Let $U_k$ denote the $k$th largest value of this sequence. Then
the extremal index $\theta$ and the probability $\pi_k$ that a cluster of the
limiting compound Poisson point process has size $k$ are given by
\begin{eqnarray} \label{extrrec}
  \theta & = & \int_0^1 P\big\{U_1=\max_{j\ge 1} W_j\le t\big\}\, dt \;=\; 1- E\min(U_1,1), \nonumber\\
  \pi_k & = & \frac{\theta_k-\theta_{k+1}}\theta \qquad \text{with}\\
  \theta_k & = & \int_0^1 P\Big\{ \sum_{j=1}^\infty 1_{(t,\infty)}(W_j)=k-1\Big\}\, dt
  \;=\;
  E\big(\min(U_{k-1},1)-\min(U_k,1)\big). \nonumber
\end{eqnarray}
Hence the extremal index is determined by the distribution of the maximum of the
geometric random walk $W_j, j\ge 1$, and the probability that a cluster of exceedances
has size $k$ is determined by the distributions of the $k$ largest ``order statistics''
of this sequence.

More recently, Gomes et al.\ (2004) (see also Gomes et al.\ (2006)) analyzed the joint
asymptotics of $k$ consecutive observations of the stationary solution of the recurrence
equation \eqref{receq}. More precisely, they proved that there exists a sequence $a_n>0$
such that
\begin{eqnarray*}
  \lim_{n\to\infty} n P\big\{ X_j>a_n x_j\; \text{ for all } 1\le j\le k\big\} & = & E
  \min_{0\le j\le k-1} (x_j^{-\kappa} W_j), \\
  \lim_{n\to\infty} n P\big\{ X_j>a_n x_j\; \text{ for some } 1\le j\le k\big\} & = & E
  \max_{0\le j\le k-1} (x_j^{-\kappa} W_j),
\end{eqnarray*}
for all $x_j>0$. Obviously, the limits on the right hand sides cannot be
expressed in terms of the extremal index $\theta$ and the cluster size
distribution $\pi_k, k\ge 1,$ for {\em all} $x_j>0$. However, this will not
even be possible in the special case that all $x_j$ are equal to some $x$, say,
so that
\begin{eqnarray*}
  \lim_{n\to\infty} n P\big\{ X_j>a_n x\; \text{ for all } 1\le j\le k\big\} & = &
  x^{-\kappa} E \min_{0\le j\le k-1} W_j \\
  &= & x^{-\kappa} \int_0^1 P\big\{ \min_{0\le j\le k-1} W_j>t\big\}\, dt, \\[0.5ex]
  \lim_{n\to\infty} n P\big\{ X_j>a_n x\; \text{ for some } 1\le j\le k\big\} & = &
  x^{-\kappa} E \max_{0\le j\le k-1}  W_j\\
  &=& x^{-\kappa} \int_0^\infty P\big\{ \max_{0\le j\le k-1} W_j>t\big\}\,
  dt,
\end{eqnarray*}
because the right hand sides depend on the distribution of the minimum and the
maximum of a {\em finite segment} of the sequence $W_j, j\ge 1,$ instead of the
distribution of the order statistics of the whole sequence.

The asymptotic variance of the Hill estimator discussed in Section 2 is another example
of a parameter which arises naturally in statistical applications and cannot be
determined from the extremal index and the cluster size distribution. Drees (2000) showed
that the Hill estimator based on the $k$ largest observations of a stationary solution of
the recurrence equation \eqref{receq} is asymptotically normal with variance
$$ \kappa^{-2}\Big( 1+2\sum_{j=1}^\infty \int_0^1 P\{W_j>t\}\, dt \Big), $$
provided that $k$ tends to $\infty$ not too fast. (An analogous result holds true for the
maximum likelihood estimator of the extreme value index, which is asymptotically normal
with a similar variance with factor $(1+1/\kappa)^2$ instead of $\kappa^{-2}$.) Here, the
parameter is determined by all marginal distributions of the geometric random walk $W_j$,
$j\ge 1$.

These examples demonstrate that in many applications the extremal index (and
often also the cluster size distribution) does not give the information one is
actually interested in. Thus there is clearly a need to put the statistical
analysis of the extremal dependence on a broader basis such that also
parameters can be treated which describe aspects of the dependence structure
different from the {\em size} of clusters of exceedances.

An interesting point of departure may be the concept of cluster functionals introduced by
Yun (2000) and developed further by Segers (2003). Roughly speaking, these are
functionals which depend only on all shortest vectors of observations containing
exceedances over a given high threshold. An asymptotic theory on estimators of
functionals of that type would be a significant step forward towards a general approach
to analyze the extremal dependence structure of stationary time series.

\section{Conclusions}

In Section 2 we compared the direct approach to the marginal tail analysis,
advocated by Rootz{\'e}n, Leadbetter and de Haan (1990), among others, with a
model-based approach where the tail of the residuals is analyzed. It was
pointed out that the perception that the latter approach is more efficient when
the model assumptions are correct is not generally justified. Moreover, it was
shown that the model-based estimators can be extremely sensitive to moderate
deviations from the models, which are difficult to detect by statistical means.
Hence, in most applications, the direct nonparametric analysis, that has proved
powerful in the classical i.i.d.\ setting in several papers by Laurens de Haan
and many others, seems also preferable for the tail analysis of serially
dependent data.

It is worth mentioning that usually the model-based approach is even more problematic if
a nonlinear time series model is assumed. For example, as we have seen in Section 3, the
marginal tail behavior of a stationary solution to the stochastic recurrence equation
\eqref{receq} does not only depend on the tail of the ``innovations'' $A_t$ (and $B_t$),
but on their whole distribution, since the extreme value index $\gamma=1/\kappa$ is
determined by the relationship $E A_1^\kappa=1$. So, in a parametric submodel, it will
not be sufficient to analyze the tail behavior of suitably defined residuals, but one has
to estimate this expectation, that depends on the center of the distribution of the
innovations and, in addition, is sensitive to deviations in the tail. Hence, to obtain a
reliable estimate, usually one has to combine some nonparametric estimate for the central
region with an extreme value estimator for the tail of the distribution of the
innovations, which makes the whole method cumbersome.

While in the marginal tail analysis sometimes too restrictive model assumptions are used,
the inference on the extremal dependence structure is often too focussed on the extremal
index (which is then estimated in quite general time series models). The nonlinear time
series \eqref{receq} is a nice example in which parameters arise in a natural way which
cannot be expressed in terms of the extremal index or the cluster size distribution. This
observation  calls for more general estimators of the extremal dependence structure.
Unfortunately, even from the most optimistic point of view, such a general theory has
just started to emerge and many challenging problems still wait for a solution.

\vspace{3ex}  {\bf Acknowledgement}: The author is grateful to J{\"u}rg H{\"u}sler and
Peng Liang for pointing out the reference Barbe and McCormick (2004).

\vspace{4ex}

{\large \bf References}

\rueck
  Ancona-Navarrete, M.A., and Tawn, J.A. (2000). A comparison of methods for estimating the extremal index.
  {\em Extremes} {\bf 3}, 5--38.

\rueck
  Barbe, Ph., McCormick, W.P. (2004).
  Tail calculus with remainder, applications to tail expansions for infinite order moving averages,
  randomly stopped sums, and related topics. {\em Extremes} {\bf  7}, 337--365.

\rueck
  Brockwell, P.J., and Davis, R.A. (2002). {\em Introduction to Time Series and
  Forecasting} (2nd ed).  Springer.

\rueck
 Datta, S., and McCormick, W.P. (1998). Inference for the tail parameters of a
linear process with heavy tail innovations. {\em Ann.\ Inst.\ Statist.\ Math.}
{\bf 50}, 337--359.

\rueck
  Davis, R.A., and Resnick, S.I. (1985). Limit theory for moving averages of random
variables with regularly varying tail probabilities. {\em Ann.\ Probab.} {\bf
13}, 179--195.


\rueck
 Dekkers, A.L.M., and de Haan, L. (1989). On the estimation of the
  extreme-value index and large quantile estimation. {\em Ann.\ Statist.}
  {\bf 17}, 1795--1832.

\rueck
 Drees, H. (2000).  Weighted approximations of tail processes for $\beta$--mixing random
variables. {\em Ann.\ Appl.\ Probab.} {\bf 10}, 1274--1301.

\rueck
 Drees, H. (2002). Tail empirical processes under mixing conditions. In: H.G.\ Dehling,
T.\ Mikosch und M.\ S{\o}rensen (eds.),
  {\em Empirical Process Techniques for Dependent Data}, 325--342, Birkh\"auser, Boston.

\rueck
 Drees, H. (2003). Extreme quantile estimation for dependent data with applications to
finance. {\em Bernoulli} {\bf 9}, 617--657.

\rueck
 Ferro, C.A.T., and Segers, J. (2003).
Inference for clusters of extreme values. {\em J.\ Roy.\ Statist.\ Soc.\ B},
{\bf 65}, 545--556.

 \rueck
 Geluk, J.L., Peng, L., and de Vries, C.G. (2000). Convolutions of
heavy-tailed random variables and applications to portfolio diversification and
${\rm MA}(1)$ time series.  {\em Adv.\ Appl.\ Probab.} {\bf 32 }, 1011--1026.

\rueck
 Geluk, J. L., and Peng, L. (2000). An adaptive optimal estimate of the tail index
for ${\rm MA}(1)$ time series.  {\em Statist.\ Probab.\ Lett.} {\bf  46},
217--227.

\rueck
  Gomes, M.I., de Haan, L., and Pestana, D. (2004). Joint exceedances of the ARCH
  process. {\em J.\ Appl.\ Probab.} {\bf 41}, 919--926.

\rueck
  Gomes, M.I., de Haan, L., and Pestana, D. (2006). Correction to: Joint exceedances of the ARCH
  process. {\em J.\ Appl.\ Probab.} {\bf 43}, 1206.

\rueck
  de Haan, L. (1990). Fighting the arch-enemy with mathematics. {\em Statist.\ Neerlandica} {\bf 44}, 45--68.

\rueck
 de Haan, L., Resnick, S.I., Rootz\'en, H., and de Vries, C. (1989).
  Extremal behaviour of solutions to a stochastic difference equation with applications to ARCH-processes.
  {\em  Stoch.\ Proc.\  Appl.} {\bf 32}, 213--224.

\rueck
  de Haan, L., and de Ronde, J. (1998). Sea and wind: multivariate extremes at work. {\em
  Extremes} {\bf 1}, 7--45.

\rueck
 de Haan, L., and Stadtm{\"u}ller, U. (1996). Generalized regular variation of second order. {\em J.\ Aust.\ Math.\ Soc.\ A} {\bf 61},
  381--395.

\rueck
 Hsing, T. (1991). On tail estimation using dependent data. {\em Ann.\ Statist.} {\bf 19},
1547--1569.

\rueck
  Hsing, T. (1993). Extremal index estimation for a weakly dependent stationary sequence. {\em
  Ann.\ Statist.} {\bf 21}, 2043--2071.

\rueck
  Hsing, T., H{\"u}sler, J., and Leadbetter, M.R. (1988). On the exceedance point process for a stationary sequence.
  {\em Probab.\ Theory
Relat. Fields} {\bf 78}, 97--112.

\rueck
  Kesten, H. (1973). Random difference equations and renewal theory for products of random matrices. {\em Acta Math.}
  {\bf 131}, 207--248.

\rueck
  Leadbetter, M.R. (1983). Extremes and local dependence in stationary sequences. {\em Probab.\ Theory
Relat. Fields} {\bf 65},
   291--306.

\rueck
  Leadbetter, M.R. (1995). On high level exceedance modeling and tail inference. {\em
J.\ Statist.\ Plann. Inference} {\bf 45}, 247--260.

\rueck
  Leadbetter, M.R., and Rootz{\'e}n, H. (1993). On central limit theory for families of strongly mixing additive random functions.
  In: Stochastic processes: a festschrift in honour of
Gopinath Kallianpur (S.\ Cambanis et al., eds.), 211--223.  Springer.

\rueck
  Ledford, A.W., and Tawn, J.A. (2003). Diagnostics for dependence within time series
extremes. {\em J.\ Royal\ Statist.\ Soc.\ B} {\bf 65}, 521-–543.

\rueck
 Ling, S., and Peng, L. (2004). Hill's estimator for the tail index of an ARMA model. {\em
J.\ Statist.\ Plann.\ Inference} {\bf 123}, 279--293.

\rueck
 Mikosch, T., and Samorodnitsky, G. (2000). The supremum of a
 negative drift random walk with dependent heavy-tailed steps. {\em Ann.\ Appl.\ Probab.} {\bf 10},
1025--1064.

\rueck
 Novak, S.Y. (2002). Inference on heavy tails from dependent data. {\em Siberian Adv.\
Math.} {\bf 12}, 73--96.

\rueck
 Pereira, T.T. (1994). Second order behaviour of domains of attraction and the bias of
generalized Pickands' estimator. In: Extreme Value Theory and Applications III
(J. Galambos, J. Lechner and E. Simiu, eds.), 165--177. NIST special
publication 866.

\rueck
 Resnick, S., and St\v{a}ric\v{a}, C. (1997). Asymptotic behavior
  of Hill's estimator for autoregressive data. {\em Comm.\ Statist.\
  Stochastic Models} {\bf 13}, 703--721.

\rueck
  Rootz\'en, H. (1995). The tail empirical process for
  stationary sequences. Unpublished manuscript, Chalmers University Gothenburg.

\rueck
  Rootz\'en, H. (2006).
Weak convergence of the tail empirical process for stationary sequences.
Submitted.

\rueck
  Rootz\'en, H., Leadbetter, M.R., and de Haan, L. (1990). Tail
  and quantile estimators for strongly mixing stationary
  processes. Report, Department
of Statistics, University of North Carolina.

\rueck
  Rootz\'en, H., Leadbetter, M.R., and de Haan, L. (1998). On the distribution of tail
  array sums for strongly mixing stationary sequences. {\em Ann.\ Appl.\ Probab.} {\bf
  8}, 868--885.

\rueck
  Segers, J. (2003). Functionals of clusters of extremes. {\em Adv.\ Appl.\ Probab.} {\bf
  35}, 1028--1045.

\rueck
  Smith, R.L., and Weissman, I. (1994). Estimating the extremal index.
  {\em J.\ Roy.\ Statist.\ Soc.\ B} {\bf 56}, 515--528.

\rueck
  Vervaat, W. (1979). On a stochastic difference equation and a representation of non-negative infinitely
divisible random variables. {\em Adv.\ Appl.\ Probab.} {\bf 11}, 750--783.

\rueck
  Weissman, I., and Novak, S.Yu. (1998). On blocks and runs estimators of the extremal index.
  {\em J.\ Statist.\ Plann.\ Inference} {\bf 66}, 281--288.

\rueck
  Yun, S. (2000). The distribution of cluster functionals of extreme events in a
  $d$th-order Markov chain. {\em J.\ Appl.\ Probab.} {\bf 37}, 29--44.

\end{document}